\documentclass[12pt]{article}
\usepackage{latexsym}
\usepackage{amssymb}
\usepackage{graphicx}
\usepackage{cite}
\usepackage{color}
\usepackage{ifpdf}

\newtheorem{Theorem}{Theorem}[part]
\newtheorem{Definition}{Definition}[part]
\newtheorem{Proposition}{Proposition}[part]

\newtheorem{Lemma}{Lemma}[part]
\newtheorem{Corollary}{Corollary}[part]

\newtheorem{Example}{Example}[part]

\def \ep{\hbox{ }\hfill$\Box$}

\def\reff#1{{\rm(\ref{#1})}}

\begin{document}
\title{The Eigenvectors of the Zero Laplacian and Signless Laplacian Eigenvalues of a Uniform Hypergraph}

\author{
Shenglong Hu \thanks{Email: Tim.Hu@connect.polyu.hk. Department of
Applied Mathematics, The Hong Kong Polytechnic University, Hung Hom,
Kowloon, Hong Kong.},\hspace{4mm} Liqun Qi \thanks{Email:
maqilq@polyu.edu.hk. Department of Applied Mathematics, The Hong
Kong Polytechnic University, Hung Hom, Kowloon, Hong Kong. This
author's work was supported by the Hong Kong Research Grant
Council (Grant No. PolyU 501909, 502510, 502111 and 501212).} }

\date{\today}
\maketitle

\begin{abstract}
In this paper, we show that the eigenvectors of the zero Laplacian
and signless Lapacian eigenvalues of a $k$-uniform hypergraph are
closely related to some configured components of that hypergraph. We
show that the components of an eigenvector of the zero Laplacian or
signless Lapacian eigenvalue have the same modulus. Moreover, under
a {\em canonical} regularization, the phases of the components of
these eigenvectors only can take some uniformly distributed values
in $\{\mbox{exp}(\frac{2j\pi}{k})\;|\;j\in [k]\}$. These
eigenvectors are divided into H-eigenvectors and N-eigenvectors.
Eigenvectors with minimal support is called {\em minimal}. The
minimal canonical H-eigenvectors characterize the even
(odd)-bipartite connected components of the hypergraph and vice
versa, and the minimal canonical N-eigenvectors characterize some
multi-partite connected components of the hypergraph and vice versa.

\vspace{2mm}
\noindent {\bf Key words:}\hspace{2mm} Tensor,
eigenvector, hypergraph, Laplacian, partition \vspace{1mm}

\noindent {\bf MSC (2010):}\hspace{2mm}
05C65; 15A18
\end{abstract}

\section{Introduction}
\setcounter{Theorem}{0} \setcounter{Proposition}{0}
\setcounter{Corollary}{0} \setcounter{Lemma}{0}
\setcounter{Definition}{0} \setcounter{Remark}{0}
\setcounter{Algorithm}{0}  \setcounter{Example}{0} \hspace{4mm} In
this paper, we study the eigenvectors of the zero Laplacian and
signless Laplacian eigenvalues of a uniform hypergraph. The zero
Laplacian eigenvalue and the zero signless Laplacian eigenvalue are
referred to respectively the zero eigenvalues of the Laplacian
tensor and the signless Laplacian tensor proposed by Qi \cite{q12a}.
It turns out that these eigenvectors are characterized by some
configured components of the underlying hypergraph and vice versa.
This work is motivated by the classic results for bipartite graphs
\cite{c97,bh11}. It is derived by the recently rapid developments on
both spectral hypergraph theory
\cite{l07,hq12a,cd12,q12a,pz12,lqy12,rp09,r09,xc12a,xc12b,xc12c,l05,hq13}
and spectral theory of tensors
\cite{cpz08,cpz11,hhlq12,hhq11,l05,l07,lqy12,hq12b,nqz09,q05,q06,q07,q12b,s12,yy10,yy11}.

The study of the Laplacian-type tensors for a uniform hypergraph becomes an active research frontier in spectral hypergraph theory recently \cite{hq12a,lqy12,xc12a,xc12c,q12a,hq13}.
Notably, Qi \cite{q12a} proposed a simple and natural definition $\mathcal
D-\mathcal A$ for the Laplacian tensor and $\mathcal D+\mathcal A$
for the signless Laplacian tensor. Here $\mathcal A=(a_{i_1\ldots
i_k})$ is the adjacency tensor of a $k$-uniform hypergraph and $\mathcal
D=(d_{i_1\ldots i_k})$ the diagonal tensor with its diagonal
elements being the degrees of the vertices. Following this, Hu and Qi proposed the normalized Laplacian tensor (or simply Laplacian) and made some explorations on it \cite{hq13}, which is the analogue of the spectral graph theory investigated by Chung
\cite{c97}.

In spectral graph theory, it is well known that the multiplicity of
the zero eigenvalue of the Laplacian matrix is equal to the number
of connected components of the graph, and the multiplicity of the
zero eigenvalue of the signless Laplacian matrix is equal to the
number of bipartite connected components of the graph \cite{bh11}.
In this paper, we investigate their analogues in spectral hypergraph
theory. It turns out that on one hand the situation is much more
complicated, and on the other hand the results are more abundant. We
show that (Please see Section 2 for the basic definitions)
\begin{itemize}
\item [(i)] Let $k$ be even and $G$ be a $k$-uniform hypergraph.
\begin{itemize}
\item (Proposition \ref{prop-s-2}) The number of minimal canonical H-eigenvectors of the zero signless Laplacian eigenvalue equals the number of odd-bipartite connected components of $G$.
\item (Proposition \ref{prop-s-5}) The number of the minimal canonical H-eigenvectors of the zero Laplacian eigenvalue equals the sum of the number of even-bipartite connected components of $G$ and the number of connected components of $G$, minus the number of singletons of $G$.
\end{itemize}

\item [(ii)] Let $k$ be odd and $G$ be a $k$-uniform hypergraph.
\begin{itemize}
\item (Proposition \ref{prop-s-3}) The number of the minimal canonical H-eigenvectors of the zero Laplacian eigenvalue equals the number of connected components of $G$.
\item (Corollary \ref{cor-1}) The number of the minimal canonical H-eigenvectors of the zero signless Laplacian eigenvalue is equal to the number of singletons of $G$.
\end{itemize}
\end{itemize}
When we turn to N-eigenvectors, we have (The definitions for the various multi-partite connected components are given in Section 6)
\begin{itemize}
\item [(i)] (Proposition \ref{prop-s-7}) Let $G=(V,E)$ be a $3$-uniform hypergraph.  Then the number of the minimal canonical conjugated N-eigenvector pairs of the zero Laplacian eigenvalue equals the number of tripartite connected components of $G$.
\item [(ii)] Let $G=(V,E)$ be a $4$-uniform hypergraph.
\begin{itemize}
\item (Proposition \ref{prop-s-8}) The number of the minimal canonical conjugated N-eigenvector pairs of the zero Laplacian eigenvalue equals the number of L-quadripartite connected components of $G$.
\item (Proposition \ref{prop-s-9}) The number of the minimal canonical conjugated N-eigenvector pairs of the zero signless Laplacian eigenvalue equals the number of sL-quadripartite connected components of $G$.
\end{itemize}
\item [(iii)] (Proposition \ref{prop-s-10}) Let $G=(V,E)$ be a $5$-uniform hypergraph. Then the number of the minimal canonical conjugated N-eigenvector pairs of the zero Laplacian eigenvalue equals the number of pentapartite connected components of $G$.
\end{itemize}
The results related with N-eigenvectors can be generalized to any
order $k\geq 6$. But it is somewhat complicated to describe the
corresponding configured components. Hence, for this part, only
$k=3,4,5$ are presented in this paper.

On top of the above results, we show in Proposition \ref{prop-s-1} that an hm-bipartite hypergraph has symmetric spectrum, and in Proposition \ref{prop-s-0} that when $k$ is even the spectrum of the Laplacian tensor and the spectrum of the signless Laplacian tensor of an hm-bipartite hypergraph equal. We also show that zero is not an eigenvalue of the signless Laplacian tensor of a connected $k$-uniform hypergraph with odd $k$ (Proposition \ref{prop-s-6}).

The rest of this paper is organized as follows.
Some definitions on eigenvalues of tensors and hypergraphs are presented in the next section. We discuss in Section 3 some spectral properties of hm-bipartite hypergraphs. Then we characterize the eigenvectors of the zero Laplacian and signless Laplacian eigenvalues of a uniform hypergraph in Section 4. In Section 5, we establish the connection between the minimal canonical H-eigenvectors of the eigenvalue zero and some configured components of the hypergraph. The discussion is extended to N-eigenvectors in Section 6.
Some final remarks are made in the last section.

\section{Preliminaries}\label{sec-p}
\setcounter{Theorem}{0} \setcounter{Proposition}{0}
\setcounter{Corollary}{0} \setcounter{Lemma}{0}
\setcounter{Definition}{0} \setcounter{Remark}{0}
\setcounter{Algorithm}{0}  \setcounter{Example}{0}
\hspace{4mm} Some preliminary definitions of eigenvalues of tensors and uniform hypergraphs are presented in this section.

\subsection{Eigenvalues of Tensors}\label{s-et}
In this subsection, some basic facts about eigenvalues and eigenvectors of tensors are reviewed. For comprehensive references, see \cite{q05,q06,q07,hhlq12} and references therein.

Let $\mathbb C$ ($\mathbb R$) be the field of complex (real) numbers and $\mathbb C^{n}$ ($\mathbb R^n$) the $n$-dimensional complex (real) space.
For integers $k\geq 3$ and $n\geq 2$, a real tensor $\mathcal T=(t_{i_1\ldots i_k})$ of order $k$ and dimension $n$ refers to a multiway array (also called hypermatrix) with entries $t_{i_1\ldots i_k}$ such that $t_{i_1\ldots i_k}\in\mathbb{R}$ for all $i_j\in[n]:=\{1,\ldots,n\}$ and $j\in[k]$. Tensors are always referred to $k$-th order real tensors in this paper, and the dimensions will be clear from the content.
Given a vector $\mathbf{x}\in \mathbb{C}^{n}$, define an $n$-dimensional vector ${\cal T}\mathbf{x}^{k-1}$ with its $i$-th element being $\sum\limits_{i_2,\ldots,i_k\in[n]}t_{ii_2\ldots i_k}x_{i_2}\cdots x_{i_k}$ for all $i\in[n]$.
Let ${\cal I}$ be the identity tensor of
appropriate dimension, e.g., $i_{i_1\ldots i_k}=1$ if and only if $i_1=\cdots=i_k\in [n]$,
and zero otherwise when the dimension is $n$.
The following definition was introduced by Qi \cite{q05}.
\begin{Definition}\label{def-00}
Let $\mathcal T$ be a $k$-th order $n$-dimensional real tensor. For
some $\lambda\in\mathbb{C}$, if polynomial system $\left(\lambda
{\cal I}-{\cal T}\right)\mathbf{x}^{k-1}=0$ has a solution
$\mathbf{x}\in\mathbb{C}^n\setminus\{0\}$, then $\lambda$ is called
an eigenvalue of the tensor ${\cal T}$ and ${\mathbf x}$ an
eigenvector of ${\cal T}$ associated with $\lambda$. If an
eigenvalue $\lambda$ has an eigenvector $\mathbf x\in\mathbb R^n$,
then $\lambda$ is called an H-eigenvalue and $\mathbf x$ an
H-eigenvector. If an eigenvector $\mathbf x\in\mathbb C^n$ cannot be
scaled to be real \footnote{We see that when $\mathbf x$ is an eigenvector, then $\alpha\mathbf x$ is still an eigenvector for all nonzero $\alpha\in\mathbb C$. An eigenvector $\mathbf x$ can be scaled to be real means that there is a nonzero $\alpha\in \mathbb C$ such that $\alpha\mathbf x\in\mathbb R^n$. In this situation, we prefer to study this $\alpha\mathbf x$, which is an H-eigenvector, other than the others in the orbit $\{\gamma\mathbf x\;|\;\gamma\in\mathbb C\setminus\{0\}\}$.}, then it is called an N-eigenvector.
\end{Definition}
It is easy to see that an H-eigenvalue is real.
However, an H-eigenvalue may still have some N-eigenvectors.

If an eigenvector $\mathbf x\in\mathbb C^n$ satisfies that some of
its components being of the maximum module one, then the eigenvector
$\mathbf x$ is {\em canonical}. In the following, unless stated
otherwise, all eigenvectors are referred to canonical eigenvectors.
This convention does not introduce any restrictions, since the
eigenvector defining equations are homogeneous. An eigenvector
$\mathbf x$ of the eigenvalue zero is called {\em minimal} if there
does not exist another eigenvector of the eigenvalue zero such that
its support is strictly contained by that of $\mathbf x$.

The {\em algebraic multiplicity} of an eigenvalue is defined as the multiplicity of this eigenvalue as a root of the characteristic polynomial $\chi_{\mathcal T}(\lambda)$. To give the definition of the characteristic polynomial, the determinant theory is needed.
For the determinant theory of a tensor, see \cite{hhlq12}.

\begin{Definition}\label{def-000}
Let $\mathcal T$ be a $k$-th order $n$-dimensional real tensor and $\lambda$ be an indeterminate variable. The determinant $\mbox{Det}(\lambda{\mathcal I}-{\cal T})$ of $\lambda{\mathcal I}-{\cal T}$, which is a
polynomial in $\mathbb{C}[\lambda]$ and denoted by $\chi_{\mathcal T}(\lambda)$, is called the {\em characteristic polynomial} of the tensor ${\cal T}$.
\end{Definition}
It is shown that the set of eigenvalues of $\mathcal T$ equals the set of roots of $\chi_{\mathcal T}(\lambda)$, see \cite[Theorem 2.3]{hhlq12}. If $\lambda$ is a root of $\chi_{\mathcal T}(\lambda)$ of multiplicity $s$, then we call $s$ the algebraic multiplicity of the eigenvalue $\lambda$. Let $c(n,k)=n(k-1)^{n-1}$. By \cite[Theorem 2.3]{hhlq12}, $\chi_{\mathcal T}(\lambda)$ is a monic polynomial of degree $c(n,k)$.
The set of all the eigenvalues of $\mathcal T$ (with algebraic multiplicities) is the {\em spectrum} of
$\mathcal T$.

Sub-tensors are involved in this paper. For more discussions on this, see \cite{hhlq12}.
\begin{Definition}\label{def-0}
Let $\mathcal T$ be a $k$-th order $n$-dimensional real tensor and $s\in [n]$. The $k$-th order $s$-dimensional tensor $\mathcal U$ with entries
$u_{i_1\ldots i_k}=t_{j_{i_1}\ldots j_{i_k}}$ for all $i_1,\ldots,i_k\in [s]$ is called the {\em sub-tensor} of ${\cal T}$
associated to the subset $S:=\{j_1,\ldots,j_s\}$. We usually denoted $\mathcal U$ as $\mathcal T(S)$.
\end{Definition}

For a subset $S\subseteq [n]$, we denoted by $|S|$ its cardinality. For $\mathbf x\in\mathbb C^n$, $\mathbf x(S)$ is
defined as an $|S|$-dimensional sub-vector of $\mathbf x$ with its
entries being $x_{i}$ for $i\in S$, and $\mbox{sup}(\mathbf x):=\{i\in[n]\;|\;x_i\neq 0\}$ is its {\em support}.

\subsection{Uniform Hypergraphs}

In this subsection, we present some preliminary concepts of uniform hypergraphs which will be used in this paper. Please refer to \cite{b73,c97,bh11,hq13,q12a} for comprehensive references.

In this paper, unless stated otherwise, a hypergraph means an undirected simple $k$-uniform hypergraph $G$ with vertex set $V$, which is labeled as $[n]=\{1,\ldots,n\}$, and edge set $E$.
By $k$-uniformity, we mean that for every edge $e\in E$, the cardinality $|e|$ of $e$ is equal to $k$. Throughout this paper, $k\geq 3$ and $n\geq k$. Moreover, since the trivial hypergraph (i.e., $E=\emptyset$) is of less interest, we consider only hypergraphs having at least one edge in this paper.

For a subset $S\subset [n]$, we denoted by $E_S$ the set of edges $\{e\in E\;|\;S\cap e\neq\emptyset\}$. For a vertex $i\in V$, we simplify $E_{\{i\}}$ as $E_i$. It is the set of edges containing the vertex $i$, i.e., $E_i:=\{e\in E\;|\;i\in e\}$. The cardinality $|E_i|$ of the set $E_i$ is defined as the {\em degree} of the vertex $i$, which is denoted by $d_i$. Then we have that $k|E|=\sum_{i\in[n]}d_i$. If $d_i=0$, then we say that the vertex $i$ is {\em isolated} or it is a {\em singleton}. Two different vertices $i$ and $j$ are {\em connected} to each other (or the pair $i$ and $j$ is connected), if there is a sequence of edges $(e_1,\ldots,e_m)$ such that $i\in e_1$, $j\in e_m$ and $e_r\cap e_{r+1}\neq\emptyset$ for all $r\in[m-1]$. A hypergraph is called {\em connected}, if every pair of different vertices of $G$ is connected. A set $S\subseteq V$ is a {\em connected component} of $G$, if every two vertices of $S$ are connected and there is no vertices in $V\setminus S$ that are connected to any vertex in $S$. For the convenience, an isolated vertex is regarded as a connected component as well. Then, it is easy to see that for every hypergraph $G$, there is a partition of $V$ as pairwise disjoint subsets $V=V_1\cup\ldots \cup V_r$ such that every $V_i$ is a connected component of $G$. Let $S\subseteq V$, the hypergraph with vertex set $S$ and edge set $\{e\in E\;|\;e\subseteq S\}$ is called the {\em sub-hypergraph} of $G$ induced by $S$. We will denoted it by $G_S$.
In the sequel, unless stated otherwise, all the notations introduced above are reserved for the specific meanings.

The following definition for the Laplacian tensor and signless Laplacian tensor was proposed by Qi \cite{q12a}.
\begin{Definition}\label{def-l}
Let $G=(V,E)$ be a $k$-uniform hypergraph. The {\em adjacency tensor} of $G$ is defined as the $k$-th order $n$-dimensional tensor $\mathcal A$ whose $(i_1 i_2\ldots i_k)$-entry is:
\begin{eqnarray*}
a_{i_1 i_2\ldots i_k}:=\left\{\begin{array}{cl}\frac{1}{(k-1)!}&if\;\{i_1,i_2\ldots,i_k\}\in E,\\0&\mbox{otherwise}.\end{array}\right.
\end{eqnarray*}
Let $\mathcal D$ be a $k$-th order $n$-dimensional diagonal tensor with its diagonal element $d_{i\ldots i}$ being $d_i$, the degree of vertex $i$, for all $i\in [n]$. Then $\mathcal D-\mathcal A$ is the {\em Laplacian tensor} of the hypergraph $G$, and $\mathcal D+\mathcal A$ is the {\em signless Laplacian tensor} of the hypergraph $G$.
\end{Definition}

Let $G=(V,E)$ be a hypergraph with connected components $V=V_1\cup\cdots\cup V_r$ for $r\geq 1$. Reorder the indices, if necessary, $\mathcal A$ can be represented by a block diagonal structure according to $V_1,\ldots, V_r$. By Definition \ref{def-00}, the spectrum of $\mathcal A$ does not change when reordering the indices. Thus, in the sequel, we assume that $\mathcal A$ is in the block diagonal structure with its $i$-th block tensor being the sub-tensor of $\mathcal A$ associated to $V_i$ for $i\in [r]$. It is easy to see that $\mathcal A(V_i)$ is the adjacency tensor of the sub-hypergraph $G_{V_i}$ for all $i\in [r]$. Similar conventions are assumed to the signless Laplacian tensor and the Laplacian tensor. By similar discussions as \cite[Lemma 3.3]{hq13}, the spectra of the adjacency tensor, the signless Laplacian tensor, and the Laplacian tensor are the unions of the spectra of its diagonal blocks respectively. For comprehensive references, please see \cite{hhlq12,hq13}.

In the following, we introduce three kinds of bipartite hypergraphs. It will be shown that the first class of uniform hypergraphs has symmetric spectra, and the latter two classes of uniform hypergraphs characterize the minimal canonical H-eigenvectors of the zero Laplacian and signless Laplacian eigenvalues.
\begin{Definition}\label{def-bi-hm}
Let $G=(V,E)$ be a $k$-uniform hypergraph. It is called {\em hm-bipartite} if either it is trivial (i.e., $E=\emptyset$) or there is a disjoint partition of the vertex set $V$ as $V=V_1\cup V_2$ such that $V_1,V_2\neq \emptyset$ and every edge in $E$ intersects $V_1$ with exactly one vertex and $V_2$ the rest $k-1$ vertices.
\end{Definition}
We use the name {\em hm-bipartite}, since a {\bf h}ead is selected
from every edge and the rest is the {\bf m}ass.

\begin{Definition}\label{def-bi-odd}
Let $k$ be even and $G=(V,E)$ be a $k$-uniform
hypergraph. It is called {\em odd-bipartite} if either it is trivial
(i.e., $E=\emptyset$) or there is a disjoint partition of the vertex
set $V$ as $V=V_1\cup V_2$ such that $V_1,V_2\neq \emptyset$ and
every edge in $E$ intersects $V_1$ with exactly an odd number of
vertices.
\end{Definition}

\begin{Definition}\label{def-bi-even}
Let $k\geq 4$ be even and $G=(V,E)$ be a $k$-uniform
hypergraph. It is called {\em even-bipartite} if either it is
trivial (i.e., $E=\emptyset$) or there is a disjoint partition of
the vertex set $V$ as $V=V_1\cup V_2$ such that $V_1,V_2\neq
\emptyset$ and every edge in $E$ intersects $V_1$ with exactly an
even number of vertices.
\end{Definition}
The idea of even-bipartite hypergraphs appeared in \cite[Corollary 6.5]{hq13}.
In Figure 1, we give preliminary examples on the three kinds of bipartite hypergraphs defined above.

\begin{figure}[htbp]
\centering
\includegraphics[width=1.6in]{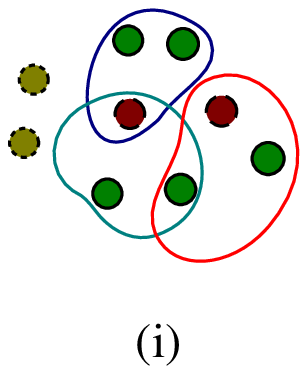}
\includegraphics[width=1.6in]{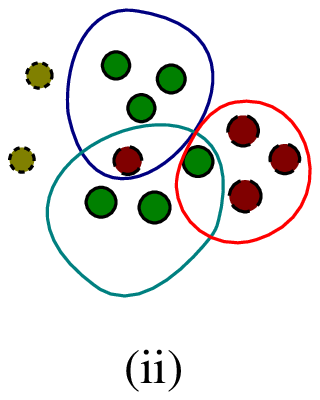}
\includegraphics[width=1.2in]{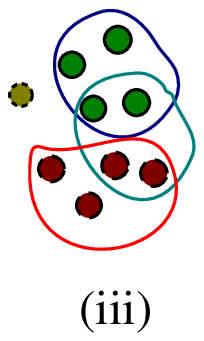}
\caption{Examples of the three kinds of bipartite hypergraphs in Definitions \ref{def-bi-hm}, \ref{def-bi-odd} and \ref{def-bi-even}. (i) is an hm-bipartite $3$-uniform hypergraph; (ii) is an odd-bipartite $4$-uniform hypergraph; and (iii) is an even-bipartite $4$-uniform hypergraph. An edge is pictured as a closed curve with the containing solid disks the vertices in that edge. Different edges are in different curves with different colors. A solid disk in dotted margin is a singleton of the hypergraph. The bipartition is clear from the different colors (also the dashed margins from the solid ones) of the disks in the connected component.}
\end{figure}

When $G$ is an ordinary graph, i.e., $k=2$, Definitions
\ref{def-bi-hm} and \ref{def-bi-odd} reduce to the classic
definition for bipartite graphs \cite{bh11}. When $k>2$, the
definition has various meaningful generalizations. In this paper, we
investigate the proposed three generalizations.

We note that not like its graph counterpart, an even (odd)-bipartite
connected component of an even-uniform hypergraph $G$ may have several
bipartitions of the same type. In this situation, we will make the
convention that this connected component contributes to the total
number of even (odd)-bipartite connected components of $G$ as the
number of the ways of the bipartitions. Likewise, for a connected
component $V_0$ of a hypergraph $G$, if it has two bipartitions of
the same type as $S_1\cup T_1=V_0$ and $S_2\cup T_2=V_0$, unless
$S_1=S_2$ or $S_1=T_2$, the two bipartitions are regarded as
different. See the following example.

\begin{Example}\label{exm-1}
{\em Let $G=(V,E)$ be a $4$-uniform hypergraph with the vertex set $V=[6]$ and the edge set $E=\{\{1,2,3,4\},\{1,3,5,6\},\{1,2,3,6\}\}$. Then, $G$ is connected and can be viewed as an even-bipartite hypergraph with a bipartition as $V_1:=\{1,2,5\}$ and $V_2:=\{3,4,6\}$ by Definition \ref{def-bi-even}. Meanwhile, $G$ is an even-bipartite hypergraph associated to the bipartitions $V_1:=\{2,3,5\}$ and $V_2:=\{1,4,6\}$; and $V_1:=\{1,3\}$ and $V_2:=\{2,4,5,6\}$. Hence the number of even-bipartite connected components of the hypergraph $G$ is three. The three even-bipartitions of the hypergraph are pictured in Figure 2.}
\end{Example}

\begin{figure}[htbp]
\centering
\includegraphics[width=1.6in]{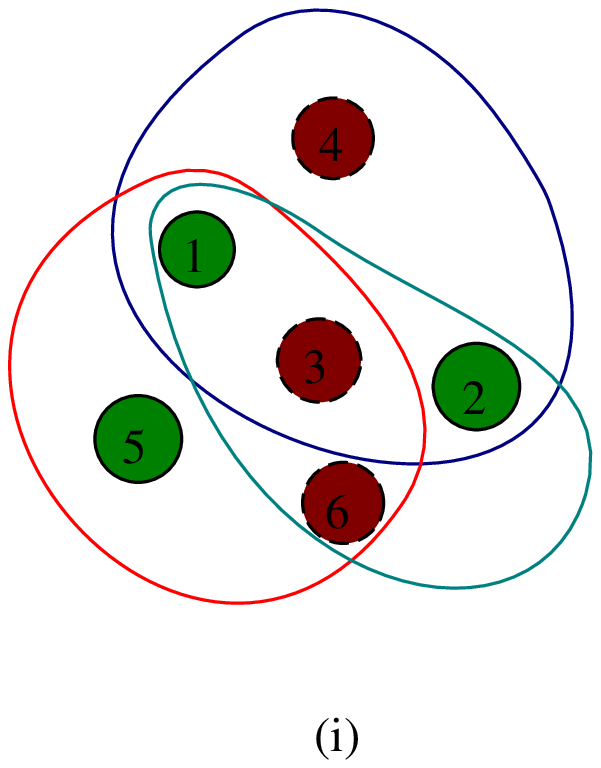}
\includegraphics[width=1.6in]{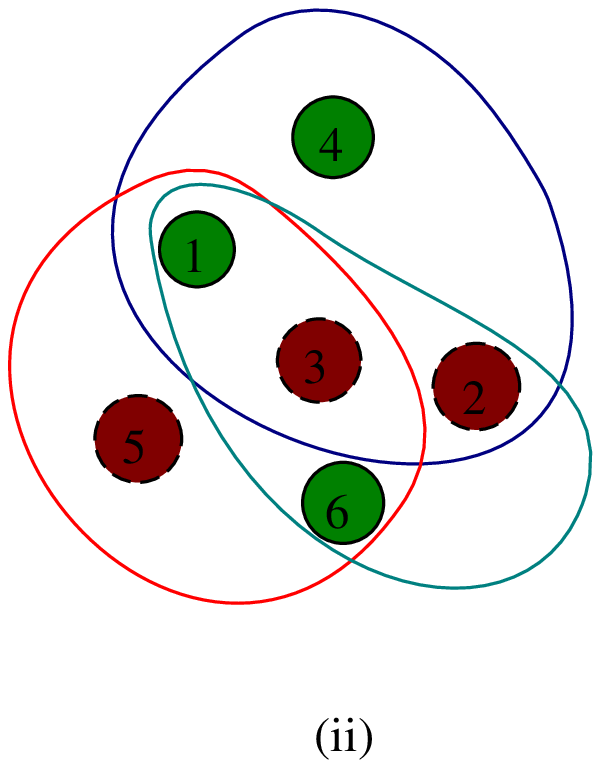}
\includegraphics[width=1.6in]{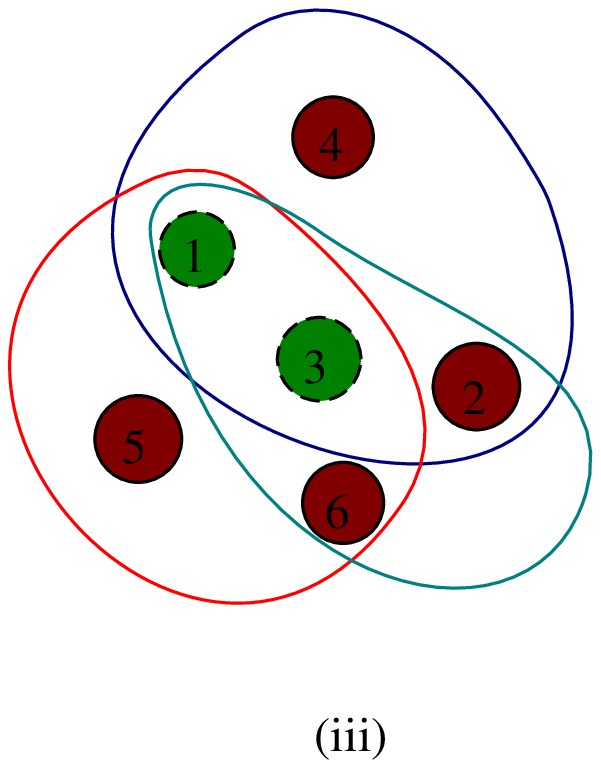}
\caption{The three even-bipartitions for the hypergraph in Example \ref{exm-1}. The legend of the pictures is similar to that of Figure 1.}
\end{figure}

\section{HM-Bipartite Hypergraphs}
\setcounter{Theorem}{0} \setcounter{Proposition}{0}
\setcounter{Corollary}{0} \setcounter{Lemma}{0}
\setcounter{Definition}{0} \setcounter{Remark}{0}
\setcounter{Algorithm}{0}  \setcounter{Example}{0} \hspace{4mm}
This section presents some basic facts about the spectra of hm-bipartite hypergraphs.

The next proposition says that the spectrum of an hm-bipartite hypergraph is symmetric. The meaning of the symmetry is clear from the statement of this proposition.

\begin{Proposition}\label{prop-s-1}
Let $G$ be a $k$-uniform hm-bipartite hypergraph. Then the spectrum of $\mathcal A$ is invariant under the multiplication of any $k$-th root of unity.
\end{Proposition}

\noindent {\bf Proof.} The case for $E=\emptyset$ is trivial. In the sequel, we assume that $E\neq\emptyset$.
Since $G$ is hm-bipartite, let $V_1$ and $V_2$ be a bipartition of $V$ such that $|e\cap V_1|=1$ for all $e\in E$.

Let $\alpha$ be any $k$-th root of unity. Suppose that $\lambda\in \mathbb C$ is an eigenvalue of $\mathcal A$ with an eigenvector being $\mathbf x\in\mathbb C^n$. Then, let $\mathbf y\in\mathbb C^n$ be a vector such that $y_i=\alpha x_i$ whenever $i\in V_1$ and $y_i=x_i$ for the others. For $i\in V_1$, we have
\begin{eqnarray*}
(\mathcal A\mathbf y^{k-1})_i=\sum_{e\in E_i}\prod_{j\in e\setminus\{i\}}y_j=\sum_{e\in E_i}\prod_{j\in e\setminus\{i\}}x_j=\lambda x_i^{k-1}=(\alpha\lambda)(\alpha x_i)^{k-1}=(\alpha \lambda) y_i^{k-1},
\end{eqnarray*}
where the second equality follows from the fact that $G$ is hm-bipartite which implies that exactly the rest vertices of every $e\in E_i$ other than $i$ belong to $V_2$, and the third from the eigenvalue equation for $(\lambda, \mathbf x)$.
For $i\in V_2$, we have
\begin{eqnarray*}
(\mathcal A\mathbf y^{k-1})_i=\sum_{e\in E_i}\prod_{j\in e\setminus\{i\}}y_j=\alpha\sum_{e\in E_i}\prod_{j\in e\setminus\{i\}}x_j=\alpha\lambda x_i^{k-1}=(\alpha\lambda)y_i^{k-1},
\end{eqnarray*}
where the second equality follows from the fact that $G$ is hm-bipartite which implies that exactly one vertex other than $i$ of every $e\in E_i$ belongs to $V_1$, and the third from the eigenvalue equation for $(\lambda, \mathbf x)$.

Hence, by Definition \ref{def-00}, $\alpha \lambda$ is an eigenvalue of $\mathcal A$. The result follows. \ep

A hypergraph is called $k$-partite, if there is a pairwise disjoint partition of $V=V_1\cup\cdots\cup V_k$ such that every edge $e\in E$ intersects $V_i$ nontrivially (i.e., $e\cap V_i\neq\emptyset$) for all $i\in [k]$. Obviously, $k$-partite hypergraphs are hm-bipartite. Thus, Proposition \ref{prop-s-1} generalizes \cite[Theorem 4.2]{cd12}.

The next proposition establishes the connection of the spectra of the signless Laplacian tensor and the Laplacian tensor for an hm-bipartite hypergraph.

\begin{Proposition}\label{prop-s-0}
Let $k$ be even and $G$ be a $k$-uniform hm-bipartite hypergraph. Then the spectrum of the Laplacian tensor and the spectrum of the signless Laplacian tensor equal.
\end{Proposition}

\noindent {\bf Proof.}
For a tensor $\mathcal T$ of order $k$ and dimension $n$, its similar transformation \footnote{This is a special matrix-tensor multiplication. Please refer to \cite{hhlq11,s12} for more properties on general matrix-tensor multiplication. } by a diagonal matrix $P$ is, denoted by $P^{-1}\cdot \mathcal T\cdot P$, defined as a $k$-th order $n$-dimensional tensor with its entries being
\begin{eqnarray*}
(P^{-1}\cdot \mathcal T\cdot P)_{i_1\ldots i_k}:=p_{i_1i_1}^{-k+1} t_{i_1\ldots i_k}p_{i_2i_2}\cdots p_{i_ki_k},\;\;\forall i_s\in [n], s\in [k].
\end{eqnarray*}
Consequently, $P^{-1}\cdot \mathcal I\cdot P=\mathcal I$. We then have
\begin{eqnarray*}
\mbox{Det}( P^{-1}\cdot (\lambda\mathcal I-\mathcal T)\cdot P)&=&\mbox{Det}(\lambda P^{-1}\cdot \mathcal I\cdot P-P^{-1}\cdot \mathcal T\cdot P)\\
&=&\mbox{Det}(\lambda \mathcal I-P^{-1}\cdot \mathcal T\cdot P)
\end{eqnarray*}
By \cite[Propositions 4.3 and 4.4]{hhlq11}, we get that
\begin{eqnarray*}
\mbox{Det}( P^{-1}\cdot (\lambda\mathcal I-\mathcal T)\cdot P)&=&[\mbox{Det}(P^{-1})]^{(k-1)^n}\mbox{Det}( \lambda\mathcal I-\mathcal T)[\mbox{Det}(P)]^{(k-1)^n}\\
&=&\mbox{Det}( \lambda\mathcal I-\mathcal T).
\end{eqnarray*}
 These two facts, together with \cite[Theorem 2.3]{hhlq12}, imply that $\mathcal T$ and $P^{-1}\cdot \mathcal T\cdot P$ have the same spectrum for any invertible diagonal matrix $P$, see also \cite{s12}.

The case for $E=\emptyset$ is trivial. Moreover, by the diagonal block structure of $\mathcal D-\mathcal A$ and $\mathcal D+\mathcal A$, we can assume that $E\neq\emptyset$ and $G$ is connected.
Since $G$ is hm-bipartite, let $V_1$ and $V_2$ be a bipartition of $V$ such that $|e\cap V_1|=1$ for all $e\in E$.
Let $P$ be a diagonal matrix with its $i$-th diagonal entry being $1$ if $i\in V_1$ and $-1$ if $i\in V_2$.
By direct computation, we have that
\begin{eqnarray*}
P^{-1}\cdot (\mathcal D-\mathcal A)\cdot P=\mathcal D+\mathcal A.
\end{eqnarray*}

Then, the conclusion follows from the preceding discussion since $P$ is invertible.  \ep

When $k$ is odd and $G$ is nontrivial, we do not have $P^{-1}\cdot
(\mathcal D-\mathcal A)\cdot P=\mathcal D+\mathcal A$. Thus, it is
unknown whether the spectrum of the Laplacian tensor and the
spectrum of the signless Laplacian tensor equal or not in this
situation. On the other hand, when $k\geq 4$ is even, at present we
are not able to prove the converse of Proposition \ref{prop-s-0},
which is well-known in spectral graph theory as: a graph is
bipartite if and only if the spectrum of the Laplacian matrix and
the spectrum of the signless Laplacian matrix equal \cite{bh11}.

There are fruitful connections between the structures of the
hypergraph and their spectra \cite{hq13,q12a}. In the sequel, the
discussions are based on the eigenvectors of the zero Laplacian and
singless Laplacian eigenvalues. To this end, we establish some basic
facts about the eigenvectors of the eigenvalue zero in the next
section.   For the Laplacian, there are similar results in
\cite[Section 6]{hq13}.

\section{Eigenvectors of the Zero Eigenvalue}
\setcounter{Theorem}{0} \setcounter{Proposition}{0}
\setcounter{Corollary}{0} \setcounter{Lemma}{0}
\setcounter{Definition}{0} \setcounter{Remark}{0}
\setcounter{Algorithm}{0}  \setcounter{Example}{0} \hspace{4mm}
The next lemma characterizes the eigenvectors of the zero Laplacian and signless Laplacian eigenvalues of a uniform hypergraph.
\begin{Lemma}\label{lem-bi-2}
Let $G$ be a $k$-uniform hypergraph and $V_i,\;i\in [s]$ be its
connected components with $s>0$. If $\mathbf x$ is an eigenvector of
the zero Laplacian or signless Laplacian eigenvalue, then $\mathbf
x(V_i)$ is an eigenvector of $(\mathcal D-\mathcal A)(V_i)$ or
$(\mathcal D+ \mathcal A)(V_i)$ corresponding to the eigenvalue zero
whenever $\mathbf x(V_i)\neq 0$. Furthermore, in this situation,
$\mbox{sup}(\mathbf x(V_i))=V_i$, and
$x_j=\gamma\mbox{exp}(\frac{2\alpha_j\pi}{k}\sqrt{-1})$ for some
nonnegative integer $\alpha_j$ for all $j\in V_i$ and some
$\gamma\in\mathbb C \setminus \{ 0 \}$.
\end{Lemma}

\noindent {\bf Proof.} The proof for the signless Laplacian tensor and that for the Laplacian tensor are similar. Hence, only the former is given.

Similar to the proof of \cite[Lemma 3.3]{hq13}, we have that for every connected component $V_i$ of $G$, $\mathbf x(V_i)$ is an eigenvector of $(\mathcal D+\mathcal A)(V_i)$ corresponding to the eigenvalue zero whenever $\mathbf x(V_i)\neq 0$.

Suppose that $\mathbf x(V_i)\neq 0$. The case for $V_i$ being a singleton is trivial. In the following, we assume that $V_i$ has more than two vertices.
We can always scale $\mathbf x(V_i)$ with some nonzero $\gamma\in\mathbb C$ such that $\frac{x_j}{\gamma}$ is positive and of the maximum module $1$ for some $j\in V_i$. Thus, without loss of generality, we assume that $\mathbf x(V_i)$ is a canonical eigenvector of $(\mathcal D+\mathcal A)(V_i)$ and $x_j=1$ for some $j\in V_i$.
Then the $j$-th eigenvalue equation is
\begin{eqnarray}\label{bi-6}
0=\left[(\mathcal D+\mathcal A) \mathbf x^{k-1}\right]_j&=&d_jx_j^{k-1}+\sum_{e\in E_j} \prod_{t\in e\setminus\{j\}}x_t=d_j+\sum_{e\in E_j} \prod_{t\in e\setminus\{j\}}x_t.
\end{eqnarray}
Since $d_j=|\{e\;|\;e\in E_j\}|$, we must have that
\begin{eqnarray*}
\prod_{t\in e\setminus\{j\}}x_t=-1,\;\forall e\in E_j.
\end{eqnarray*}
This implies that
\begin{eqnarray}\label{bi-7}
\prod_{t\in e}x_t=-1,\;\forall e\in E_j.
\end{eqnarray}
Since the maximum module is $1$, $x_t=\mbox{exp}(\theta_t\sqrt{-1})$ for some $\theta_t\in [0,2\pi]$ for all $t\in e$ with $e\in E_j$.
For another vertex $s$ which shares an edge with $j$, we have
\begin{eqnarray*}
0=\left[(\mathcal D+\mathcal A) \mathbf x^{k-1}\right]_s&=&d_sx_s^{k-1}+\sum_{e\in E_s} \prod_{t\in e\setminus\{s\}}x_t.
\end{eqnarray*}
Similarly, we have
\begin{eqnarray*}
x_s^{k-1}=-\prod_{t\in e\setminus\{s\}}x_t,\;\forall e\in E_s.
\end{eqnarray*}
Thus,
\begin{eqnarray*}
x_s^k=-\prod_{t\in e}x_t,\;\forall e\in E_s.
\end{eqnarray*}
The fact that $s$ and $j$ share one edge, together with \reff{bi-7},
implies that $x_s^k=1$. Similarly, we have that
\begin{eqnarray*}
x_s^k=1,\;\forall s\in e,\;e\in E_j.
\end{eqnarray*}
As $G_{V_i}$ is connected, by induction, we can show that $x_s^k=1$ for all $s\in V_i$. Consequently, $\theta_t=\frac{2\alpha_t}{k}\pi$ for some integers $\alpha_t$ for all $t\in V_i$. \ep

With Lemma \ref{lem-bi-2}, parallel results as those in \cite[Section 6]{hq13} can be established for the signless Laplacian tensor and the Laplacian tensor. Especially, we have the following result.

\begin{Theorem}\label{thm-h-2}
Let $G$ be a $k$-uniform connected hypergraph.
\begin{itemize}
\item [(i)] A nonzero vector $\mathbf x$ is an eigenvector of the Laplacian tensor $\mathcal D-\mathcal A$ corresponding to the zero eigenvalue if and only if there exist nonzero $\gamma\in\mathbb C$ and integers $\alpha_i$ such that $x_i=\gamma\exp(\frac{2\alpha_i\pi}{k}\sqrt{-1})$ for $i\in[n]$, and
\begin{eqnarray}\label{bi-8}
\sum_{j\in e}\alpha_j=\sigma_e k,\;\forall e\in E,
\end{eqnarray}
for some integer $\sigma_e$ with $e\in E$.
\item [(ii)] A nonzero vector $\mathbf x$ is an eigenvector of the signless Laplacian tensor $\mathcal D+\mathcal A$ corresponding to the zero eigenvalue if and only if there exist nonzero $\gamma\in\mathbb C$ and integers $\alpha_i$ such that $x_i=\gamma\exp(\frac{2\alpha_i\pi}{k}\sqrt{-1})$ for $i\in[n]$, and
\begin{eqnarray}\label{bi-9}
\sum_{j\in e}\alpha_j=\sigma_e k+\frac{k}{2},\;\forall e\in E,
\end{eqnarray}
for some integer $\sigma_e$ with $e\in E$.
\end{itemize}
\end{Theorem}

An immediate consequence of Theorem \ref{thm-h-2} is that the singless Laplacian tensor of a $k$-uniform connected hypergraph with odd $k$ does not have zero eigenvalue, since $\frac{k}{2}$ is a fraction and \reff{bi-9} can never be fulfilled in this situation.
\begin{Proposition}\label{prop-s-6}
Let $k$ be odd and $G$ be a $k$-uniform connected hypergraph. Then zero is not an eigenvalue of the signless Laplacian tensor.
\end{Proposition}

When we restrict the discussion of Lemma \ref{lem-bi-2} to H-eigenvectors, we get the following corollary. We state it explicitly here for the convenience of the next section's reference.
\begin{Corollary}\label{lem-bi-1}
Let $G$ be a $k$-uniform hypergraph and $V_i,\;i\in [s]$ be its connected components with $s>0$. If $\mathbf x$ is an H-eigenvector of the zero Laplacian or signless Laplacian eigenvalue, then $\mathbf x(V_i)$ is an H-eigenvector of $(\mathcal D-\mathcal A)(V_i)$ or $(\mathcal D+\mathcal A)(V_i)$ corresponding to the eigenvalue zero whenever $\mathbf x(V_i)\neq 0$. Furthermore, in this situation, $\mbox{sup}(\mathbf x(V_i))=V_i$, and up to a real scalar multiplication $x_j=\pm 1$ for $j\in V_i$.
\end{Corollary}
By Corollary \ref{lem-bi-1}, we get that a minimal canonical H-eigenvector $\mathbf x\in\mathbb R^n$ of the zero Laplacian or signless Laplacian eigenvalue must be of this form: $x_i$ is $1$, $-1$ or $0$ for all $i\in [n]$.

\section{H-Eigenvectors}
\setcounter{Theorem}{0} \setcounter{Proposition}{0}
\setcounter{Corollary}{0} \setcounter{Lemma}{0}
\setcounter{Definition}{0} \setcounter{Remark}{0}
\setcounter{Algorithm}{0}  \setcounter{Example}{0} \hspace{4mm} In this section, we establish the connection between the H-eigenvectors of the zero Laplacian and signless Laplacian eigenvalues and the even (odd)-connected components of the underlying hypergraph.

We note that canonical eigenvectors are considered in this section. Then, when we do the number count, we always consider a minimal canonical H-eigenvector $\mathbf x$ and its reciprocal $-\mathbf x$ as the same.

The next proposition, together with Corollary \ref{lem-bi-1}, generalizes \cite[Theorem 1.3.9]{bh11} which says that the multiplicity of the eigenvalue zero of the signless Laplacian matrix of a graph is equal to the number of bipartite connected components of this graph.

\begin{Proposition}\label{prop-s-2}
Let $k$ be even and $G$ be a $k$-uniform hypergraph. Then the number of minimal canonical H-eigenvectors of the zero signless Laplacian eigenvalue equals the number of odd-bipartite connected components of $G$.
\end{Proposition}

\noindent {\bf Proof.} Suppose that $V_1\subseteq V$ is an odd-bipartite connected component of $G$.
If $V_1$ is a singleton, then $1$ is a minimal canonical H-eigenvector of $(\mathcal D+\mathcal A)(V_1)=0$ by definition. In the following, we assume that $G_{V_1}$ has at least one edge.
Let $V_1=S\cup T$ be an odd-bipartition of the sub-hypergraph $G_{V_1}$ such that every edge of $G_{V_1}$ intersects with $S$ an odd number of vertices. Then $S, T\neq \emptyset$, since $k$ is even. Let $\mathbf y\in \mathbb R^{|V_1|}$ be a vector such that $y_i=1$ whenever $i\in S$ and $y_i=-1$ whenever $i\in T$. Then, for $i\in S$,
\begin{eqnarray*}
\left[(\mathcal D+\mathcal A) \mathbf y^{k-1}\right]_i&=&d_iy_i^{k-1}+\sum_{e\in E_i} \prod_{j\in e\setminus\{i\}}y_j\\
&=&d_i-d_i\\
&=&0.
\end{eqnarray*}
Here the second equality follows from the fact that for every $e\in E_i$ exactly the number of $y_j=-1$ is odd.

Next, for $i\in T$,
\begin{eqnarray*}
\left[(\mathcal D+\mathcal A) \mathbf y^{k-1}\right]_i&=&d_iy_i^{k-1}+\sum_{e\in E_i} \prod_{j\in e\setminus\{i\}}y_j\\
&=&-d_i+d_i\\
&=&0.
\end{eqnarray*}
Here the second equality follows from the fact that for every $e\in E_i$ here exactly the number of $y_j=-1$ other than $y_i$ is even.

Thus, for every odd-bipartite connected component of $G$, we can associate it a canonical H-eigenvector corresponding to the eigenvalue zero. Since $d_i=|\{e\;|\;e\in E_i\}|$, it is easy to see that this vector $\mathbf y$ is a minimal canonical H-eigenvector. Otherwise, suppose that $\mathbf x$ with $\mbox{sup}(\mathbf x)\subset\mbox{sup}(\mathbf y)$ is a canonical H-eigenvector of $(\mathcal D+\mathcal A)(V_1)$ corresponding to the eigenvalue zero.
Since $V_1$ is a nontrivial connected component of $G$ and $\mathbf x\neq 0$, we must have a $j\in V_1$ such that $x_j\neq 0$ and there is an edge containing both $j$ and $s$ with $x_s=0$.
The $j$-th eigenvalue equation is
\begin{eqnarray*}
0=\left[(\mathcal D+\mathcal A) \mathbf y^{k-1}\right]_j&=&d_jx_j^{k-1}+\sum_{e\in E_j} \prod_{t\in e\setminus\{j\}}x_t.
\end{eqnarray*}
We have $|\prod_{t\in e\setminus\{j\}}x_t|\leq 1$ and $|x_j|=1$. Since $x_s=0$ and there is an edge containing both $s$ and $j$, we have $|\sum_{e\in E_j} \prod_{t\in e\setminus\{j\}}x_t|\leq |\{e\;|\;e\in E_j\}|-1<|\{e\;|\;e\in E_j\}|=d_j$. Thus, a contradiction to the eigenvalue equation. Hence, $\mathbf y$ is minimal.
Obviously, if $S_1\cup T_1=V_1$ and $S_2\cup T_2=V_1$ are two different odd-bipartitions of the connected component $V_1$, then the constructed minimal canonical H-eigenvectors are different.
So the number of odd-bipartite connected components of $G$ is not greater than the number of minimal canonical H-eigenvectors of the zero signless Laplacian eigenvalue.

Conversely, suppose that $\mathbf x\in \mathbb R^n$ is a minimal canonical H-eigenvector corresponding to the eigenvalue zero, then $\mbox{sup}(\mathbf x)$ is a connected component of $G$ by Corollary \ref{lem-bi-1}. Denote by this connected component of $G$ as $V_0$. If $V_0$ is a singleton, then it is an odd-bipartite connected component by Definition \ref{def-bi-odd}. In the following, we assume that $V_0$ has more than one vertex.

For all $j\in V_0$,
\begin{eqnarray}\label{bi-2}
0=\left[(\mathcal D+\mathcal A) \mathbf y^{k-1}\right]_j&=&d_jx_j^{k-1}+\sum_{e\in E_j} \prod_{s\in e\setminus\{j\}}x_s.
\end{eqnarray}
Let $S\cup T=V_0$ be a bipartition of $V_0$ such that $x_s>0$ whenever $s\in S$ and $x_s<0$ whenever $s\in T$. Since $\mathbf x$ is canonical and $|V_0|>1$, we must have $S\neq\emptyset$. This, together with \reff{bi-2}, implies that $T\neq \emptyset$. From \reff{bi-2}, we see that for every edge $e\in E_j$ with $j\in S$, $|e\cap T|$ must be an odd number; and for every edge $e\in E_j$ with $j\in T$, $|e\cap T|$ must be an odd number as well. Then $V_0$ is an odd-bipartite component of $G$. Hence, every minimal canonical H-eigenvector corresponding to the eigenvalue zero determines an odd-bipartite connected component of $G$. Obviously, the odd-bipartite connected components determined by a minimal canonical H-eigenvector $\mathbf x$ and its reciprocal $-\mathbf x$ are the same.

Combining the above results, the conclusion follows. \ep

The next proposition is an analogue of Proposition \ref{prop-s-2} for the Laplacian tensor.
\begin{Proposition}\label{prop-s-5}
Let $k$ be even and $G$ be a $k$-uniform hypergraph. Then the number of the minimal canonical H-eigenvectors of the zero Laplacian eigenvalue equals the sum of the number of even-bipartite connected components of $G$ and the number of connected components of $G$, minus the number of singletons of $G$.
\end{Proposition}

\noindent {\bf Proof.} By similar proof of that for Proposition \ref{prop-s-2}, we see that if $V_0$ is an even-bipartite connected component of $G$, then, we can construct a minimal canonical H-eigenvector of $\mathcal D-\mathcal A$. By the proof and the definition of even-bipartite hypergraph, we see that this minimal canonical H-eigenvector contains negative components whenever this connected component is not a singleton. Meanwhile, by Definition \ref{def-00} and Lemma \ref{lem-bi-1}, for every connected component of $G$, the vector of all ones is a minimal canonical H-eigenvector corresponding to the eigenvalue zero. The two systems of minimal canonical H-eigenvectors overlap at the singletons of $G$. Thus, the number of the minimal canonical H-eigenvectors of the zero Laplacian eigenvalue is not smaller than the sum of the number of even-bipartite connected components of $G$ and the number of connected components of $G$, minus the number of singletons of $G$.

Conversely, suppose that $\mathbf x$ is a minimal canonical H-eigenvectors of the zero Laplacian eigenvalue. Then, by Lemma \ref{lem-bi-1} and similar proof for that in Proposition \ref{prop-s-2}, $\mbox{sup}(\mathbf x)$ is a connected component of $G$. If $\mathbf x$ contains both positive and negative components, then $\mbox{sup}(\mathbf x)$ must be an even-bipartite connected component. When $\mathbf x$ contains only positive components, we can only conclude that $\mbox{sup}(\mathbf x)$ is a connected component. The overlap also happens when $\mbox{sup}(\mathbf x)$ is a singleton. Thus, we can see that the number of the minimal canonical H-eigenvectors of the zero Laplacian eigenvalue is not greater than the sum of the number of even-bipartite connected components of $G$ and the number of connected components of $G$, minus the number of singletons of $G$.

Combining these two conclusions, the result follows. \ep

The next proposition is an analogue of Proposition \ref{prop-s-5} for $k$ being odd.
\begin{Proposition}\label{prop-s-3}
Let $k$ be odd and $G$ be a $k$-uniform hypergraph. Then the number of the minimal canonical H-eigenvectors of the zero Laplacian eigenvalue equals the number of connected components of $G$.
\end{Proposition}

\noindent {\bf Proof.} By a similar proof of that in Proposition \ref{prop-s-5}, we can see that a minimal canonical H-eigenvector of the eigenvalue zero can only have positive components, since $k$ is odd. Then, the result follows. \ep

The next proposition follows from Proposition \ref{prop-s-6}. We highlight it here to make the statements in this section more complete.
\begin{Proposition}\label{prop-s-4}
Let $k$ be odd and $G$ be a $k$-uniform connected hypergraph. Then the signless Laplacian tensor has no zero H-eigenvalue.
\end{Proposition}

The next corollary is an analogue of Proposition \ref{prop-s-3} for the signless Laplacian tensor, which follows from Proposition \ref{prop-s-4}.
\begin{Corollary}\label{cor-1}
Let $k$ be odd and $G$ be a $k$-uniform hypergraph. Then the number of the minimal canonical H-eigenvectors of the zero signless Laplacian eigenvalue is equal to the number of singletons of $G$.
\end{Corollary}

\section{N-Eigenvectors}
\setcounter{Theorem}{0} \setcounter{Proposition}{0}
\setcounter{Corollary}{0} \setcounter{Lemma}{0}
\setcounter{Definition}{0} \setcounter{Remark}{0}
\setcounter{Algorithm}{0}  \setcounter{Example}{0} \hspace{4mm}
In this section, we establish the connection between the minimal canonical N-eigenvectors of the zero Laplacian and signless Laplacian eigenvalues and some multi-partite connected components of the underlying hypergraph.

We note that canonical N-eigenvectors are considered in this section. Then, when we do the number count, we always consider a minimal canonical eigenvector $\mathbf x$ and its multiplication $\mbox{exp}(\theta\sqrt{-1})\mathbf x$ by any $\theta\in [0,2\pi]$ as the same. Actually, by Lemma \ref{lem-bi-2}, only $\theta=\frac{2\alpha\pi}{k}$ for integers $\alpha$ are involved.

\subsection{$3$-Uniform Hypergraphs}
In this subsection, $3$-uniform hypergraphs are discussed.

\begin{Definition}\label{def-3-1}
Let $G=(V,E)$ be a $3$-uniform hypergraph. If there exists a partition of $V=V_1\cup V_2\cup V_3$ such that $V_1,V_2,V_3\neq \emptyset$, and for every edge $e\in E$, either $e\subseteq V_i$ for some $i\in [3]$ or $e$ intersects $V_i$ nontrivially for all $i\in [3]$, then $G$ is called {\em tripartite}.
\end{Definition}
We prefer to use {\em tripartite} to distinguish it from the concept of {\em 3-partite}. Note that if a connected component $V_0$ of $G$ is tripartite, then $|V_0|\geq3$. As that for bipartite hypergraphs introduced in Section 2, for a connected tripartite hypergraph $G$ which has two tripartitions as $V=V_1\cup V_2\cup V_3=S_1\cup S_2\cup S_3$, unless (with a possible renumbering of the subscripts of $S_i,\;i\in [3]$) $S_j=V_j,\;j\in [3]$, the two tripartitions are regarded as different. Please see Figure 3 as an example. In that figure, the partition (i) and the partition (iii) are considered as different.

The next lemma helps to establish the connection.
\begin{Lemma}\label{lem-bi-3}
Let $G=(V,E)$ be a $k$-uniform hypergraph. Then $\mathbf x\in\mathbb C^n$ is a canonical N-eigenvector of the zero Laplacian or signless Laplacian eigenvalue if and only if the conjugate $\mathbf x^H$ of $\mathbf x$ is a canonical N-eigenvector of the zero Laplacian or signless Laplacian eigenvalue.
\end{Lemma}

\noindent {\bf Proof.} Since the Laplacian tensor and the signless Laplacian tensor of a uniform hypergraph are real. The conclusion follows from the definition of eigenvalues immediately. \ep

The next corollary is a direct consequence of Lemma \ref{lem-bi-3}.
\begin{Corollary}\label{cor-2}
Let $G=(V,E)$ be a $3$-uniform hypergraph. If $\mathbf x\in\mathbb C^n$ is a canonical N-eigenvector of the zero Laplacian eigenvalue, then
\begin{eqnarray*}
x_i^H&=&1,\;\forall i\in \{j\;|\;x_j=1\},\\
x_i^H&=&\mbox{exp}(\frac{4\pi}{3}\sqrt{-1}),\;\forall i\in \{j\;|\;x_j=\mbox{exp}(\frac{2\pi}{3}\sqrt{-1})\},\;\mbox{and}\\
x_i^H&=&\mbox{exp}(\frac{2\pi}{3}\sqrt{-1}),\;\forall i\in \{j\;|\;x_j=\mbox{exp}(\frac{4\pi}{3}\sqrt{-1})\}.
\end{eqnarray*}
Consequently, the tripartition of $[n]$ determined by $\mathbf x$ and that by $\mathbf x^H$ are the same.
\end{Corollary}
Corollary \ref{cor-2} can be generalized to general $k\geq 4$ routinely.

With the above lemmas and the definitions, we are in the position to present the following result.
\begin{Proposition}\label{prop-s-7}
Let $G=(V,E)$ be a $3$-uniform hypergraph. Then the number of the minimal canonical conjugated N-eigenvector pairs of the zero Laplacian eigenvalue equals the number of tripartite connected components of $G$.
\end{Proposition}

\noindent {\bf Proof.} If $V_0$ is a tripartite connected component of $G$, then let $V_0=R\cup S\cup T$ be a tripartition of it.
Let vector $\mathbf y\in\mathbb C^{|V_0|}$ such that $y_i=1$ whenever $i\in R$; $y_i=\mbox{exp}(\frac{2\pi}{3}\sqrt{-1})$ whenever $i\in S$; and $y_i=\mbox{exp}(\frac{4\pi}{3}\sqrt{-1})$ whenever $i\in T$. By Definition \ref{def-3-1} and Theorem \ref{thm-h-2} (i), we see that $\mathbf y$ is an N-eigenvector of $(\mathcal D-\mathcal A)(V_0)$ corresponding to the eigenvalue zero. By Lemma \ref{lem-bi-2}, we see that $\mathbf y$ embedded trivially in $\mathbb C^n$ is a minimal canonical N-eigenvector. Likewise, by a similar proof of that in Proposition \ref{prop-s-2}, we can conclude this result. By Lemma \ref{lem-bi-3}, the conjugate of $\mathbf y$ is also a minimal canonical N-eigenvector. Moreover, the minimal canonical N-eigenvector constructed in this way from the tripartition $V_0=R\cup S\cup T$ only can be either $\mathbf y$ or $\mathbf y^H$.
Hence the number of the minimal canonical conjugated N-eigenvector pairs of the zero Laplacian eigenvalue is not less than the number of tripartite connected components of $G$.

Conversely, we assume that $\mathbf x\in\mathbb C^n$ is a minimal canonical N-eigenvector of the zero Laplacian eigenvalue. By Lemma \ref{lem-bi-2}, $V_0:=\mbox{sup}(\mathbf x)$ is a connected component of $G$. Since $\mathbf x$ is an N-eigenvector, $V_0$ is a nontrivial connected component of $G$, i.e., $|V_0|>1$. This, together with Lemma \ref{lem-bi-2}, implies that we can get a tripartition of $V_0$ as $V_0=R\cup S\cup T$ such that $R$ consists of the indices such that the corresponding components of $\mathbf x$ being one; $S$ consists of the indices such that the corresponding components of $\mathbf x$ being $\mbox{exp}(\frac{2\pi}{3}\sqrt{-1})$; and $T$ consists of the indices such that the corresponding components of $\mathbf x$ being $\mbox{exp}(\frac{4\pi}{3}\sqrt{-1})$.

It is easy to see that $R\neq \emptyset$, and $S\cup T\neq \emptyset$ since $\mathbf x$ is an N-eigenvector. The fact that $\mathbf x$ is an eigenvector corresponding to the eigenvalue zero, together with Theorem \ref{thm-h-2} (i), implies that both $S$ and $T$ are nonempty. Otherwise,
the equations in \reff{bi-8} cannot be fulfilled. By \reff{bi-8} again, we must have that for every edge $e\in E(G_{V_0})$ ($i.e., E_{V_0}$) either $e\subseteq R$, $S$ or $T$; or $e$ intersects $R$, $S$ and $T$ nontrivially. Consequently, $V_0$ is a tripartite connected component of $G$. We also see from Definition \ref{def-3-1}, Lemma \ref{lem-bi-3} and Corollary \ref{cor-2} that the tripartite connected component of $G$ determined by $\mathbf x^H$ is the same as that of $\mathbf x$.
Hence, the number of the minimal canonical conjugated N-eigenvector pairs of the zero Laplacian eigenvalue is not greater than the number of tripartite connected components of $G$.

Combining the above two conclusions, the result follows. \ep

By Proposition \ref{prop-s-6}, the signless Laplacian tensor of a $3$-uniform connected hypergraph does not have the eigenvalue zero. Hence, it is wondering that whether all the configured components of $3$-uniform connected hypergraphs are characterized by the eigenvectors of the Laplacian tensor or not.

In the following, we give an example to illustrate the result obtained.
\begin{Example}\label{exm-2}
{\em Let $G=(V,E)$ be a $3$-uniform hypergraph with $V=[7]$ and
\begin{eqnarray*}
E=\{\{1,2,3\},\{3,4,5\},\{5,6,7\}\}.
\end{eqnarray*}
By Definition \ref{def-3-1}, there are three tripartitions of $G$ as
\begin{eqnarray*}
\{1\},\{2\},\{3,4,5,6,7\};\;\{1,2,3\},\{4\},\{5,6,7\};\;\mbox{and}\; \{1,2,3,4,5\},\{6\},\{7\}.
\end{eqnarray*}
These tripartitions are pictured in Figure 3. By Proposition \ref{prop-s-7}, there are three minimal canonical conjugated N-eigenvector pairs of the zero Laplacian eigenvalue. By Proposition \ref{prop-s-3}, there is only one minimal H-eigenvector of the zero Laplacian eigenvalue, i.e., the vector of all ones, since $G$ is connected.}
\end{Example}

\begin{figure}[htbp]
\centering
\includegraphics[width=1.5in]{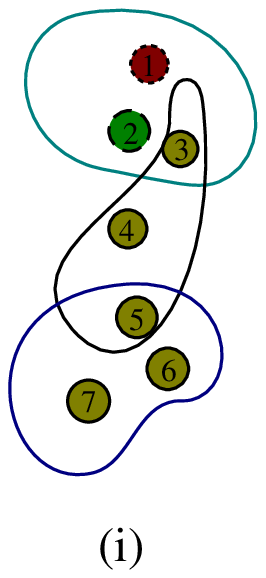}
\includegraphics[width=1.5in]{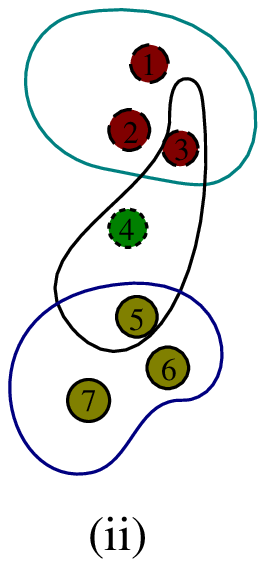}
\includegraphics[width=1.5in]{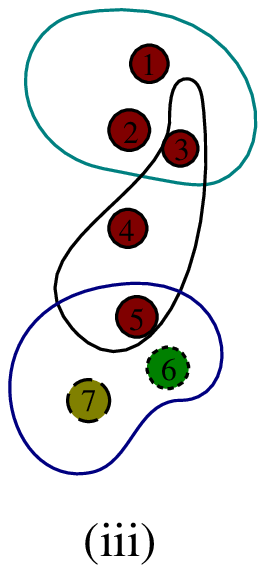}
\caption{The three tripartitions of the $3$-uniform hypergraph in Example \ref{exm-2}. The tripartitions are clear from the groups of disks in different colors (also with dotted, dashed and solid margins).}
\end{figure}

We note that, geometrically, the tripartitions in (i) and (iii) of Figure 3 are essentially the same. The numbering which is used to distinguish the tripartitions (i) and (iii) in this paper is manually added.
Then, there is a gap between Proposition \ref{prop-s-7} and the complete characterization of the intrinsic tripartite connected components of a uniform hypergraph.

We try to generalize the result in this subsection to $k$-uniform hypergraphs with $k\geq 4$. However, the definitions for multi-partite connected components of $k$-uniform hypergraphs with bigger $k$ are somewhat complicated. Then only $k=4$ and $k=5$ are presented in the next two subsections respectively. More neat statements are expected in the future.

\subsection{$4$-Uniform Hypergraphs}

In this subsection, $4$-uniform hypergraphs are discussed.

\begin{Definition}\label{def-4-1}
Let $G=(V,E)$ be a $4$-uniform hypergraph. If there exists a partition of $V=V_1\cup V_2\cup V_3\cup V_4$ such that at least two of $V_i, i\in [4]$ are nonempty, and for every edge $e\in E$, either $e\subseteq V_i$ for some $i\in [4]$ or one of the following situations happens:
\begin{itemize}
\item [(i)] $|e\cap V_1|=2$ and $|e\cap V_3|=2$;
\item [(ii)] $|e\cap V_2|=2$ and $|e\cap V_4|=2$;
\item [(iii)] $|e\cap V_1|=2$, $|e\cap V_2|=1$, and $|e\cap V_4|=1$;
\item [(iv)] $|e\cap V_3|=2$, $|e\cap V_2|=1$, and $|e\cap V_4|=1$,
\end{itemize}
then $G$ is called {\em L-quadripartite}.
\end{Definition}
Similarly, we prefer to use {\em L-quadripartite} to distinguish it from the concept of {\em 4-partite}. Here the prefix ``L" is for the Laplacian tensor.

\begin{Definition}\label{def-4-1}
Let $G=(V,E)$ be a $4$-uniform hypergraph. If there exists a partition of $V=V_1\cup V_2\cup V_3\cup V_4$ such that at least two of $V_i, i\in [4]$ are nonempty, and for every edge $e\in E$, either $e\subseteq V_i$ for some $i\in [4]$ or one of the following situations happens:
\begin{itemize}
\item [(i)] $|e\cap V_1|=3$ and $|e\cap V_3|=1$;
\item [(ii)] $|e\cap V_3|=3$ and $|e\cap V_1|=1$;
\item [(iii)] $|e\cap V_2|=3$ and $|e\cap V_4|=1$;
\item [(iv)] $|e\cap V_4|=3$ and $|e\cap V_2|=1$;
\item [(v)] $|e\cap V_1|=2$ and $|e\cap V_2|=2$;
\item [(vi)] $|e\cap V_2|=2$ and $|e\cap V_3|=2$;
\item [(vii)] $|e\cap V_3|=2$ and $|e\cap V_4|=2$;
\item [(viii)] $e$ intersects $V_i$ nontrivially for all $i\in [4]$,
\end{itemize}
then $G$ is called {\em sL-quadripartite}.
\end{Definition}
Here the prefix ``sL" is for the signless Laplacian tensor. For a connected L(sL)-quadripartite hypergraph $G$ which has two L(sL)-quadripartitions as $V=V_1\cup \cdots\cup V_4=S_1\cup\cdots\cup S_4$, unless (with a possible renumbering of the subscripts of $S_i,\;i\in [4]$) $S_j=V_j,\;j\in [4]$, the two L(sL)-quadripartitions are regarded as different.

With similar proofs as that for Proposition \ref{prop-s-7}, we can get the following two propositions. However, we give the key points of the proof for the next proposition as an illustration.
\begin{Proposition}\label{prop-s-8}
Let $G=(V,E)$ be a $4$-uniform hypergraph. Then the number of the minimal canonical conjugated N-eigenvector pairs of the zero Laplacian eigenvalue equals the number of L-quadripartite connected components of $G$.
\end{Proposition}

\noindent {\bf Proof.} Let $V_0$ be an L-quadripartite connected component  of $G$ with the L-quadripartition being $S_1\cup S_2\cup S_3\cup S_4$. Without loss of generality, we can assume that $S_1\neq \emptyset$.
Then we can associated it a minimal canonical N-eigenvector $\mathbf x\in\mathbb C^n$ through
\begin{eqnarray*}
x_i:=\mbox{exp}(\frac{2(j-1)\pi}{4}\sqrt{-1}),\;\mbox{whenever}\;i\in S_j,\;j\in [4].
\end{eqnarray*}
By Definition \ref{def-4-1}, we can conclude that the other minimal canonical N-eigenvectors $\mathbf y$ with $y_i=1$ for $i\in S_1$ can only be $\mathbf y=\mathbf x^H$. Consequently, an L-quadripartite connected component of $G$ determines a minimal canonical conjugated N-eigenvector pair. The converse is obvious. Then the result follows with a similar proof as that for Proposition \ref{prop-s-7}.  \ep

\begin{Proposition}\label{prop-s-9}
Let $G=(V,E)$ be a $4$-uniform hypergraph. Then the number of the minimal canonical conjugated N-eigenvector pairs of the zero signless Laplacian eigenvalue equals the number of sL-quadripartite connected components of $G$.
\end{Proposition}

\subsection{$5$-Uniform Hypergraphs}

In this subsection, $5$-uniform hypergraphs are discussed.

\begin{Definition}\label{def-5-1}
Let $G=(V,E)$ be a $5$-uniform hypergraph. If there exists a partition of $V=V_1\cup\cdots\cup V_5$ such that at least three of $V_i, i\in [5]$ are nonempty, and for every edge $e\in E$, either $e\subseteq V_i$ for some $i\in [5]$ or one of the following situations happens:
\begin{itemize}
\item [(i)] $|e\cap V_2|=2$, $|e\cap V_5|=2$, and $|e\cap V_1|=1$;
\item [(ii)] $|e\cap V_3|=2$, $|e\cap V_4|=2$, and $|e\cap V_1|=1$;
\item [(iii)] $|e\cap V_1|=3$, $|e\cap V_2|=1$, and $|e\cap V_5|=1$;
\item [(iv)]  $|e\cap V_1|=3$, $|e\cap V_3|=1$, and $|e\cap V_4|=1$;
\item [(v)]  $|e\cap V_2|=3$, $|e\cap V_4|=1$, and $|e\cap V_5|=1$;
\item [(vi)] $|e\cap V_3|=3$, $|e\cap V_1|=1$, and $|e\cap V_4|=1$;
\item [(vii)] $|e\cap V_4|=3$, $|e\cap V_1|=1$, and $|e\cap V_2|=1$;
\item [(viii)] $|e\cap V_5|=3$, $|e\cap V_2|=1$, and $|e\cap V_3|=1$;
\item [(ix)] $e$ intersects $V_i$ nontrivially for all $i\in [5]$,
\end{itemize}
then $G$ is called {\em pentapartite}.
\end{Definition}
For a connected pentapartite hypergraph $G$ which has two pentapartitions as $V=V_1\cup \cdots\cup V_5=S_1\cup\cdots\cup S_5$, unless (with a possible renumbering of the subscripts of $S_i,\;i\in [5]$) $S_j=V_j,\;j\in [5]$, the two pentapartitions are regarded as different. With similar proof as that for Propositions \ref{prop-s-7} and \ref{prop-s-8}, we can get the next proposition.

\begin{Proposition}\label{prop-s-10}
Let $G=(V,E)$ be a $5$-uniform hypergraph. Then the number of the minimal canonical conjugated N-eigenvector pairs of the zero Laplacian eigenvalue equals the number of pentapartite connected components of $G$.
\end{Proposition}

\section{Final Remarks}
\setcounter{Theorem}{0} \setcounter{Proposition}{0}
\setcounter{Corollary}{0} \setcounter{Lemma}{0}
\setcounter{Definition}{0} \setcounter{Remark}{0}
\setcounter{Algorithm}{0}  \setcounter{Example}{0} \hspace{4mm} In
this paper, the relations of the eigenvectors of the zero Laplacian
and signless Laplacian eigenvalues of a uniform hypergraph, with
some configured components of that hypergraph, are discussed. It is
different from the recent work \cite{q12a,hq13} which mainly
concentrates on the discussions of H$^+$-eigenvalues of the
Laplacian tensor, the signless Laplacian tensor and the Laplacian.
H-eigenvectors and, more importantly, N-eigenvectors are discovered
to be applicable in spectral hypergraph theory. It is shown that the
H-eigenvectors and N-eigenvectors of the zero Laplacian and signless
Laplacian eigenvalues can characterize some intrinsic structures of
the underlying hypergraph. More work on the other eigenvalues of
these Laplacian-type tensors are expected in the future. Especially,
a general statement on the connection between the N-eigenvectors and
the configured components of a hypergraph is our next topic.

{\bf Acknowledgement.} The authors are very grateful to Prof. Jia-Yu Shao for valuable comments.


\bibliographystyle{model6-names}

\end{document}